\documentclass[12pt]{amsart}

\usepackage{environ, xcolor}
\NewEnviron{commentA}{}
\newcommand\Aon{\RenewEnviron{commentA}{\color{purple}\BODY}}

\Aon{} %Change this to Aoff for professional version, Aon to see extra comments

\usepackage{amsmath, amssymb}
\usepackage{array}
\usepackage[frame,cmtip,arrow,matrix,line,graph,curve]{xy}
\usepackage{graphpap, color, paralist, pstricks}
\usepackage[mathscr]{eucal}
\usepackage[pdftex]{graphicx}
\usepackage[pdftex,colorlinks,backref=page,citecolor=blue]{hyperref}
\usepackage{cleveref}
\usepackage{tikz, verbatim}

\newtheorem{theorem}{Theorem}[section]
\newtheorem{proposition}[theorem]{Proposition}

\newtheorem{lemma}[theorem]{Lemma}

\newtheorem{conjecture}[theorem]{Conjecture}

\theoremstyle{definition}
\newtheorem{definition}[theorem]{Definition}

\newtheorem{remark}[theorem]{Remark}

\usepackage[margin=1in]{geometry} 
\usepackage{amsmath,amssymb,mathtools, mathrsfs,esint,tikz-cd,verbatim}
\newcommand{\R}{{\mathbb R}}

\newcommand{\Z}{{\mathbb Z}}
\newcommand{\Q}{{\mathbb Q}}

\newcommand{\C}{{\mathbb C}}

\newcommand{\on}[1]{\operatorname{#1}}

\newcommand{\Spec}{{\on{Spec}}}


\title{Minimal log discrepancies of hypersurface mirrors}
\author{Louis Esser}

\address{Department of Mathematics, Princeton University, Fine Hall, Washington Road, Princeton, NJ 08544-1000, USA}
\email{esserl@math.princeton.edu}

\begin{document}
	\maketitle

 \begin{abstract}
For certain quasismooth Calabi--Yau hypersurfaces in weighted projective space, the Berglund-H\"{u}bsch-Krawitz (BHK) mirror symmetry construction gives a concrete description of the mirror.  We prove that the minimal log discrepancy of the quotient of such a hypersurface by its toric automorphism group is closely related to the weights and degree of the BHK mirror.  As an application, we exhibit klt Calabi--Yau varieties with the smallest known minimal log discrepancy.  We conjecture that these examples are optimal in every dimension.
 \end{abstract}

\section{Introduction}

We say that a normal projective variety $X$ is \textit{Calabi--Yau} if the canonical class $K_X$ is $\Q$-linearly equivalent to zero.  Similarly, a pair $(X,D)$ is Calabi--Yau if $K_X + D \sim_{\Q} 0$.  Many theorems and conjectures in algebraic geometry assert that the numerical invariants of Calabi--Yau pairs with relatively mild singularities and a fixed dimension are bounded in various ways.  For instance, the volume of an ample Weil divisor on a klt Calabi--Yau pair of dimension $n$ with coefficients in a DCC set $I$ is known to have a positive lower bound depending only on $n$ and $I$, by a result of Birkar \cite[Corollary 1.4]{Birkar}.  The \textit{index conjecture} predicts that the \textit{index} of $X$, the smallest positive integer $m$ such that $m(K_X+D) \sim 0$, is uniformly bounded for the same class of pairs.  

A series of papers by Totaro, Wang, and the author \cite{ETW, ETWindex} laid out examples of Calabi--Yau varieties and pairs that are particularly extreme with respect to the two invariants mentioned, volume and index.  Several of these examples are in fact conjecturally optimal with respect to the invariant of interest.  Surprisingly, the conjecturally optimal examples for volume and index are related by mirror symmetry \cite[Remark 3.8]{ETWindex}.  In this paper, we describe another connection between mirror symmetry and the birational geometry of Calabi--Yau pairs, this time relating to the minimal log discrepancy.

When a pair $(X,D)$ has relatively mild singularities, the minimal log discrepancy (or mld for short) is a way of quantifying how singular it is.  A pair is Kawamata log terminal (klt) if the minimal log discrepancy is positive.  It follows from results of Hacon, M\textsuperscript{c}Kernan, and Xu that there is a positive lower bound on the mld among all klt Calabi--Yau pairs of a fixed dimension $n$ with coefficients in a fixed DCC set $I$ (see \Cref{mldbound}).

The framework of mirror symmetry suggests that there is a mirror dual associated to a Calabi--Yau variety.  In certain special cases, an explicit construction for the mirror is known.  The example we'll focus on in this paper is Berglund-H\"{u}bsch-Krawitz (BHK) mirror symmetry, which describes how to find the mirror for certain orbifolds which are quotients of Calabi--Yau hypersurfaces in weighted projective space.  

Our main result shows that the mld of the ``maximal toric quotient" of any Calabi--Yau weighted projective hypersurface to which BHK mirror symmetry applies has a remarkably simple description in terms of the mirror.  This gives a new connection between mirror symmetry and invariants from birational geometry.

\begin{theorem}[\Cref{mainmldthm}]
    Let $V_d \subset \mathbb{P}_{\mathbb{C}}(a_0,\ldots,a_{n+1})$ be a well-formed quasismooth Calabi--Yau hypersurface defined by a polynomial equation of degree $d$ with the same number of monomials as variables such that its matrix of exponents is invertible (also called a Delsarte polynomial). Suppose that $\mathrm{Aut}_T(V)$ is the toric automorphism group of $V$. Let $d^{\mathsf{T}}$ and $a_0^{\mathsf{T}},\ldots,a_{n+1}^{\mathsf{T}}$ be the mirror degree and mirror weights of $V$, respectively.  Then, the minimal log discrepancy of the quotient pair $V/\mathrm{Aut}_T(V)$ is
\begin{equation}
\label{mainmldequality}
\mathrm{mld}(V/\mathrm{Aut}_T(V)) = \frac{\min\{a_0^{\mathsf{T}},\ldots,a_{n+1}^{\mathsf{T}}\}}{d^{\mathsf{T}}}.
\end{equation}
\end{theorem}

Using this result, we compute the mld for some special examples.  For instance, we write down examples of klt Calabi--Yau varieties (rather than pairs) with the smallest known mld.  Some of the properties of these examples are explicated in a separate paper \cite{ET}.  Their mld decreases doubly exponentially with dimension, and we expect that they achieve the smallest mld of any klt Calabi--Yau variety in each dimension (\Cref{kltvarconj}).  This conjecture is supported by low-dimensional evidence.  These results complement known examples of klt Calabi--Yau pairs with standard coefficients due to Jihao Liu \cite[Remark 2.6]{Liu} that have similar asymptotics.  We also show how the properties of these latter examples can be deduced as a special case of \Cref{mainmldthm}.

The key idea of the proof of \Cref{mainmldthm} will be to view quotients of hypersurfaces as above as hypersurfaces in fake weighted projective stacks, a special case of toric Deligne-Mumford stacks.  We then use the geometric description of the mld for toric singularities to verify \eqref{mainmldequality}.  The outline of the paper is as follows.  In \Cref{background_section}, we explain the necessary background on hypersurfaces in fake weighted projective stacks and their singularities; these results generalize more familiar ones for hypersurfaces in weighted projective space, but to the author's knowledge have not appeared in the literature.  The end of \Cref{background_section} also summarizes the construction of the Berglund-H\"{u}bsch-Krawitz mirror.  \Cref{main_theorem_section} is devoted to the proof of \Cref{mainmldthm}, while \Cref{applications_section} applies this result by constructing examples of klt Calabi--Yau pairs and varieties with the smallest known minimal log discrepancy.

\noindent{\it Acknowledgements. } The author was partially
supported by NSF grant DMS-2054553. Thanks to Burt Totaro for useful conversations, and to
the referees for many helpful suggestions and comments.

 \section{Notation and Background}
 \label{background_section}

 Throughout the paper, we'll work over the complex numbers $\C$.

 \subsection{Minimal log discrepancies and toric stacks}

The \textit{minimal log discrepancy} (or mld) of a pair $(X,D)$ is a numerical measure of its singularities.  

For us, a pair $(X,D)$ consists of a normal projective variety $X$ and an effective $\Q$-divisor $D$ with the property that $K_X + D$ is a $\Q$-Cartier divisor.  Then, for a proper birational morphism $\mu: X' \rightarrow X$, where $X'$ is again normal, and any irreducible divisor $E \subset X'$, we may define the \textit{log discrepancy} of $E$ (with respect to $\mu$) as follows:
$$a_E(X,D) \coloneqq \mathrm{ord}_E(K_{X'} + E - \mu^* (K_X + D)).$$
The log discrepancy $a_E(X,D)$ only depends on the valuation defined by $E$ on the function field of $X$ and not on the particular birational model $\mu: X' \rightarrow X$ for $E$.  The \textit{center} $c_X(E)$ of $E$ in $X$ is the image $\mu(E) \subset X$, which again depends only on the valuation.  For any point $x$ of the scheme $X$, the \textit{minimal log discrepancy} of the pair $(X,D)$ at the point $x$ is defined as 
$$\mathrm{mld}_x(X,D) \coloneqq \inf \{a_E(X,D): c_X(E) = \bar{x} \}.$$
The (global) \textit{minimal log discrepancy} of $(X,D)$ is
$$\mathrm{mld}(X,D) \coloneqq \inf_{x \in X}\mathrm{mld}_x(X,D),$$
where the infimum is taken over all points $x$ of the scheme $X$.  Whenever the pair $(X,D)$ is \textit{log canonical} (lc), that is, whenever $\mathrm{mld}(X,D)$ is nonnegative, the global mld can be computed using a single log resolution of $(X,D)$ and hence is rational \cite[Definition 7.1]{Kollarsing}.  A pair $(X,D)$ is \textit{Kawamata log terminal} (klt) if $\mathrm{mld}(X,D) > 0$.

Recall that $(X,D)$ is \textit{Calabi--Yau} if $K_X + D \sim_{\Q} 0$. For klt Calabi--Yau pairs $(X,D)$ with $X$ of a fixed dimension $n$ and $D$ with coefficients belonging to a fixed set $I$ satisfying the descending chain condition, there is a positive lower bound on the minimal log discrepancy:

\begin{proposition}
\label{mldbound}
Let $n$ be a positive integer and $I$ a DCC set.  Then there is a positive number $\epsilon = \epsilon(n,I)$ such that every klt Calabi--Yau pair $(X,D)$ with dimension $n$ and coefficients in $I$ has mld at least $\epsilon$.
\end{proposition}

This follows from work of Hacon-M\textsuperscript{c}Kernan-Xu \cite{HMX}; see \cite[Proposition 2.1]{ETWindex} for a proof, which uses ideas from \cite[Lemma 3.13]{CDHJS}.

In this paper, we'll focus primarily on klt Calabi--Yau pairs with \textit{standard coefficients}, meaning that $I = \{0\} \cup \{1-\frac{1}{b}: b \in \Z_+\}$.  These pairs arise naturally as certain quotients of varieties by finite groups: more precisely, if $Y$ is a normal projective variety with $\Q$-Cartier canonical class and an action by a finite group $G$, then the variety $X \coloneqq Y/G$ is naturally equipped with a divisor $D$ such that $(X,D)$ is a pair with standard coefficients.  This $D$ has the property that $K_Y = \pi^*(K_X + D)$, where $\pi: Y \rightarrow X$ is the quotient morphism.  The divisor is determined from the $G$-action in the sense that $D$ has coefficient $1 - \frac{1}{b}$ on the image of a prime divisor in $Y$ for which the subgroup of $G$ acting as the identity on the divisor has order $b$.

Throughout the paper, we'll need the explicit description of the minimal log discrepancy of toric pairs.  A \textit{toric pair} is a pair $(X,D)$ with $X$ a normal $\Q$-factorial toric variety and $D$ a torus-invariant $\Q$-divisor.  When describing the fans of toric varieties and stacks, we'll use the notation $\mathrm{cone}\{v_1,\ldots,v_k\}$ to mean the cone generated by vectors $v_1, \ldots, v_k \in \mathbb{R}^N$, i.e., the set of points which are a nonnegative linear combination of these vectors.

It will be geometrically convenient to phrase our results in the language of toric Deligne-Mumford stacks \cite[Section 3]{BCS}.  In the same way that the datum of a simplicial fan corresponds to a normal, $\Q$-factorial toric variety, the datum of a stacky fan corresponds to a toric Deligne-Mumford stack.  A \textit{stacky fan} $\mathbf{\Sigma}$ consists of a triple $(N,\Sigma,\beta)$. Here $N$ is a finitely generated abelian group, $\Sigma$ is a full-dimensional and strictly convex rational simplicial fan in $\bar{N} \coloneqq N \otimes_{\Z} \R$ with $r$ rays $\rho_1,\ldots,\rho_r$, and $\beta$ is a collection of $r$ elements $\{\beta_1,\ldots,\beta_r\}$ of $N$ such that the image $\bar{\beta}_i$ of $\beta_i$ in $\bar{N}$ spans $\rho_i$.  The list $\{\beta_1,\ldots,\beta_r\}$ is the same as a map $\beta: \Z^r \rightarrow N$ with a finite cokernel, where the $i$th standard basis element of $\Z^r$ maps to $\beta_i$. 

The Deligne-Mumford stack associated to $\mathbf{\Sigma}$ will be denoted $\mathcal{Y}(\mathbf{\Sigma})$.  We'll only explain the construction of this stack when $N$ is torsion-free, which corresponds to the stack having trivial generic stabilizer. This assumption will hold for the rest of the paper. We can then view $N \cong \Z^d$ as a lattice inside $N_{\mathbb{R}}$, and for simplicity we'll use the notation $\beta_i$ rather than $\bar{\beta}_i$ for the images of elements of $N$ distinguished by $\beta$. To construct the stack in that case, we take a quotient of an open subset of $\mathbb{A}^r$, with coordinates $z_1,\ldots,z_r$, as follows.  Let $J_{\Sigma}$ be the monomial ideal generated by the set of products $\prod_{i: \rho_i \notin \sigma} z_i$, where $\sigma$ ranges over all cones in $\Sigma$.  Denote by $\beta^{\star}: M \rightarrow (\Z^r)^{\star}$ the dual to $\beta$ (where $M$ is the dual lattice to $N$) and consider the following exact sequence:
$$0 \rightarrow M \rightarrow (\Z^r)^{\star} \rightarrow S(\mathbf{\Sigma}) \rightarrow 0.$$
Here $S(\mathbf{\Sigma})$ is a finitely generated abelian group which is the cokernel of the map $\beta^{\star}$. Then $G(\mathbf{\Sigma}) \coloneqq \Spec(\C[S(\mathbf{\Sigma})])$ carries the structure of an algebraic group, canonically embedded in $(\C^*)^r = \Spec(\C[(\Z^r)^{\star}])$ via the surjection $\C[(\Z^r)^{\star}] \rightarrow \C[S(\mathbf{\Sigma})]$.  The group $G$ also has an action on $\mathbb{A}^r$ inherited from the diagonal action of $(\C^*)^r$ on this space.  Then $\mathcal{Y}(\mathbf{\Sigma})$ is then defined as the quotient stack $[(\mathbb{A}^r \setminus V(J_{\Sigma}))/G]$.  This is in fact a Deligne-Mumford stack \cite[Proposition 3.2]{BCS}.

By \cite[Proposition 3.7]{BCS}, the coarse moduli space of the stack $\mathcal{Y}(\mathbf{\Sigma})$ is the ordinary toric variety $Y \coloneqq Y(\Sigma)$ associated to the fan $\Sigma$ in $N$.  Thus, there is a natural toric pair structure $(Y,D)$ associated to the stack $\mathcal{Y}$, where the coefficient of the torus-invariant divisor corresponding to $\rho_i$ in $D$ is $1 - \frac{1}{m_i}$; here, $m_i$ denotes the order of the subgroup which acts as the identity on the divisor under the $G(\mathbf{\Sigma})$-action.  The \textit{minimal log discrepancy} of the stack $\mathcal{Y}$ is by definition the mld of the associated pair $(Y,D)$.

Since the fan $\Sigma$ is rational simplicial, there exists a unique piecewise linear function $\psi: |\Sigma| \rightarrow \R$ which has value $1$ on each $\beta_i$ and is linear on each cone $\sigma \in \Sigma$.  We'll call this the \textit{log discrepancy function} of the toric stack $\mathcal{Y}$.  The lemma below justifies this terminology.  It is the analog for toric stacks of the geometric description of the mld of toric singularities due to A. Borisov \cite[section 2]{ABorisov}.

\begin{lemma}
Let $\mathcal{Y}(\mathbf{\Sigma})$ be a toric Deligne-Mumford stack with trivial generic stabilizer and let $\psi$ be the log discrepancy function of $\mathcal{Y}$.  Then the minimal log discrepancy of $\mathcal{Y}$ is
$$\on{mld}(\mathcal{Y}) = \min\{\psi|_{(N \cap |\Sigma|) \setminus \{0\}} \}.$$
\end{lemma}

\begin{proof}
    Since the mld is a local invariant, it suffices to consider the affine case where there is a single cone in the fan $\Sigma$.  If $d$ is the rank of $N$, then the cone has $r = d$ rays since it is simplicial; these rays are spanned by $\beta_1,\ldots,\beta_d$, respectively.  The image of the map $\beta: \Z^d \rightarrow N$ is a finite index sublattice of $N$.  We'll use $e_1, \ldots, e_d$ for the standard basis of $\Z^d$, so $\beta(e_i) = \beta_i$.   In this case, $V(J_\Sigma) = \emptyset$ and the stack $\mathcal{Y}$ is the quotient of $\mathbb{A}^d$ by the finite group $G = \Spec(\C[S(\mathbf{\Sigma})])$.  This is equivalent to the usual toric description of a (possibly ill-formed) affine abelian quotient singularity.
    
    Let $p_i$ be the primitive lattice point on the ray spanned by $\beta_i$.  We claim that the order of the subgroup of $G$ acting as the identity on $\{z_i = 0\} \subset \mathbb{A}^d$ is precisely the ratio $\beta_i/p_i$. Indeed, the subscheme $G$ of $(\C^*)^d = \C[(\Z^d)^{\star}]$ is defined by the ideal $\C[M] \subset \C[(\Z^d)^{\star}]$. The subgroup acting as the identity on $\{z_i = 0\}$ in $(\C^*)^d$ is defined by the ideal $\C[H_i]$, where $H_i \subset (\Z^d)^{\star}$ is the subgroup generated by the dual standard basis vectors $f_1,\ldots,\hat{f}_i,\ldots,f_d$ in $(\Z^d)^{\star}$.  The intersection of $G$ with this subgroup is then the spectrum of the group algebra of 
    \begin{equation}
    \label{cokernel}
         \mathrm{coker}(M \xrightarrow{\beta^{\star}} (\Z^d)^{\star} \rightarrow (\Z^d)^{\star}/H_i).
    \end{equation}
    The quotient $(\Z^d)^{\star}/H_i$ is the dual of the sublattice $\Z \cdot e_i \subset \Z^d \subset N$.  Therefore, the image of the composition in \eqref{cokernel} can be thought of as the collection of linear maps $\Z \cdot e_i \rightarrow \Z$ which extend to the ambient lattice $N$. The cokernel therefore has order equal to $\beta_i/p_i$, so the spectrum of its group algebra is a finite discrete algebraic group of that same order.

    Since the log discrepancy function $\psi$ is $1$ at $\beta_i$, it has value $p_i/\beta_i$ at the primitive lattice point $p_i$.  In addition, $1 - p_i/\beta_i$ is the coefficient of the image of the divisor $\{z_i = 0\}$ in the pair $(Y,D)$ associated to $\mathcal{Y}$.  The usual log discrepancy function of a toric pair is defined by precisely these conditions, that is, by linearity on every cone of the fan and a value of $1 - \mathrm{coeff}_{D_i} D$ on the primitive lattice point spanning the ray associated to each torus-invariant divisor $D_i$ \cite[Section 1]{Ambro}.  Therefore, this calculation confirms that the function $\psi$ deserves the name ``log discrepancy function": for any primitive lattice point $e \in N \cap |\Sigma|$, the prime divisor $E_e$ over $X$ corresponding to the barycentric subdivision of $\Sigma$ with center $e$ has log discrepancy $a_{E_e}(Y,D) = \psi(e)$.  The minimum value of $\psi$ on $(N \cap |\Sigma|) \setminus \{0\}$ is therefore the mld of $(Y,D)$, as required.
\end{proof}

\subsection{Hypersurfaces in Fake Weighted Projective Stacks}
\label{fwps}

In the main theorem, we'll be interested in computing the mld of certain pairs which are quotients of hypersurfaces in weighted projective space by finite groups.  Since these pairs are quotients of varieties with canonical singularities by finite groups, they are klt, and therefore have positive mld. It will be useful for us to view these same quotients as hypersurfaces inside of what we call fake weighted projective stacks.

A \textit{fake weighted projective stack} is a toric Deligne-Mumford stack $\mathcal{Y}(\mathbf{\Sigma}) = (N,\Sigma,\beta)$ for which $N \cong \Z^{n+1}$ and the fan $\Sigma$, which is complete and generated by $r = n+2$ rays $\rho_0, \ldots, \rho_{n+1}$, defines a $\Q$-factorial projective toric variety of Picard number $1$.  The $\beta_i$ are not required to be primitive lattice points on these rays, so this construction is slightly more general than the usual definition of fake weighted projective \textit{spaces}, which in turn generalize ordinary weighted projective space $\mathbb{P}(a_0,\ldots,a_{n+1})$.  This subsection will generalize some of the theory of weighted projective hypersurfaces and their singularities (see \cite{Iano-Fletcher} or \cite[Section 2]{ETW}) to this more general setting.

The top-dimensional cones of the fan $\Sigma$ will be denoted $\sigma_i$, $i = 0,\ldots, n+1$.  Here $\sigma_i \coloneqq \mathrm{cone}\{\beta_0,\ldots,\hat{\beta}_i,\ldots,\beta_{n+1}\}$.  We'll use the coordinates $x_0,\ldots,x_{n+1}$ for affine space $\mathbb{A}^{n+2}$.  It follows that $J_{\Sigma}$ has generators $x_0,\ldots,x_{n+1}$, so that $\mathcal{Y} = [(\mathbb{A}^{n+2} \setminus \{0\})/G]$. In this case, $G = \Spec(\C[S(\mathbf{\Sigma})])$ is a one-dimensional algebraic group of multiplicative type (or equivalently, the product of $\mathbb{G}_{\on{m}}$ with a finite abelian group).

There is a unique collection of positive integers $a_0,\ldots,a_{n+1}$ for which $a_0 \beta_0 + \cdots + a_{n+1} \beta_{n+1} = 0$ and $\gcd(a_0,\ldots,a_{n+1}) = 1$.  These are the \textit{weights} of the fake weighted projective stack $\mathcal{Y}$.  The toric variety corresponding to $\Sigma$ in the lattice generated by $\beta_0,\ldots,\beta_{n+1}$ is simply the usual weighted projective space $\mathbb{P}(a_0,\ldots,a_{n+1})$, so $\mathcal{Y}$ is naturally the quotient of this space by a finite group.  To refer to a point of the coarse moduli space $Y$ of $\mathcal{Y}$, we sometimes use homogeneous coordinates $(x_0: \cdots : x_{n+1})$.

The action of $G$ on $\mathbb{A}^{n+2}$ also gives an action of $G$ on the polynomial ring $\C[x_0,\ldots,x_{n+1}]$, which in turn gives a grading on this ring by the group of characters $S \coloneqq S(\mathbf{\Sigma})$ of $G$ (this is the same $S(\mathbf{\Sigma})$ as above, namely the cokernel of the dual of $\beta$).  Suppose that $\chi: G \rightarrow \C^*$ is a character of $G$.  An element $f \in \C[x_0,\ldots,x_{n+1}]$ is called \textit{homogeneous of degree} $\chi$ if $g \cdot f = \chi(g) f$ for all $g \in G$.  For example, when the stack $\mathcal{Y}$ equals $\mathbb{P}(a_0,\ldots,a_{n+1})$, the group of characters is simply $\Z$ and ``homogeneous of degree $d \in \Z$" means ``homogeneous of weighted degree $d$" in the usual sense.  Since the action of $G$ is diagonal, every monomial is homogeneous of some degree.  We'll write $\theta_0,\ldots,\theta_{n+1}$ for the characters of $x_0,\ldots,x_{n+1}$, respectively.  The datum consisting of the group $G$ and these \textit{coordinate characters} $\theta_0,\ldots,\theta_{n+1}$ also determines the fake weighted projective stack $\mathcal{Y}$.

For any subset $I \subset \{0,\ldots,n+1\}$ of size $k$, there is an associated toric stratum $W_I \subset \mathbb{A}^{n+2} \setminus \{0\}$ where precisely the coordinates in $I$ are nonvanishing.  We have $W_I \cong (\C^*)^{k}$.  We'll use $U_I$ to denote the image of this stratum in $Y$, which is also the set where the homogeneous coordinates indexed by $I$ are nonzero; then $U_I \cong (\C^*)^{k-1}$.

In this language, it's straightforward to identify the quotient singularities of a fake weighted projective stack $\mathcal{Y}$ in terms of the coordinate characters $\theta_0,\ldots,\theta_{n+1}$. To write down these singularities, we use the following notation.  Let $H$ be a finite group and $\chi_1,\ldots,\chi_l$ be characters of this group.  A pair $(V,D_V)$ has a \textit{quotient singularity of type} $\frac{1}{H}(\chi_1,\ldots,\chi_l)$ at a closed point $p \in V$ if there is an \'{e}tale neighborhood of $p$ which is isomorphic to the quotient pair $\mathbb{A}^l/H$, where $H$ acts diagonally by characters $\chi_1,\ldots,\chi_l$.  We'll often abuse notation slightly and say that the smooth Deligne-Mumford stack $\mathcal{V}$ with trivial generic stabilizer and coarse moduli pair $(V,D_V)$ has a quotient singularity of the given type at $p \in V$.  

Now, we can precisely describe the singularities of $\mathcal{Y}$ (cf. \cite[Proposition 2.3]{ETW}):

\begin{proposition}
    Let $\mathcal{Y} = [(\mathbb{A}^{n+2} \setminus \{0\})/G]$ be a fake weighted projective stack with coordinate characters $\theta_0,\ldots,\theta_{n+1}$ and coarse moduli pair $(Y,D_Y)$. Let $I \subset \{0,\ldots,n+1\}$ be a subset of size $|I| = k$ and $G_I$ be the intersection $\bigcap_{i \in I} \ker(\theta_i)$.  Then at any point $p$ of the toric stratum $U_I \subset Y$ where exactly the coordinates indexed by $I$ are nonvanishing, $\mathcal{Y}$ has quotient singularity
    $$\frac{1}{G_I}(\theta_i|_{G_I}: i \notin I) \times \mathbb{A}^{k-1}.$$
\end{proposition}

\begin{proof}
Let $q$ be a preimage of $p$ in $\mathbb{A}^{n+2}$.  Precisely the coordinates of $q$ indexed by $I$ are nonvanishing, so the stabilizer of $q$ is the subgroup $G_I$ of $G$ on which all coordinate characters $\theta_i, i \in I$ are trivial.  The stratum $U_I \subset Y$ is isomorphic to $(\C^*)^{k-1}$ and has the same singularities at every point.  The action of $G_I$ on the coordinates not in $I$ gives the quotient singularity shown.
\end{proof}

When $f$ is homogeneous of degree some character $\chi$, it defines a hypersurface $\mathcal{X}$ in $\mathcal{Y}$; indeed, let $C_f \coloneqq \{f = 0\} \subset \mathbb{A}^{n+2}$ be the \textit{affine cone} of $f$ and $C_f^* \coloneqq C_f \setminus \{0\}$ the \textit{punctured} affine cone.  Then $C_f^*$ is $G$-invariant, so there is a stack $\mathcal{X} = [C_f^*/G]$.  When $C_f^*$ is smooth, we say that $\mathcal{X}$ is a \textit{quasismooth} hypersurface in $\mathcal{Y}$.  When $\mathcal{X}$ is quasismooth, it is a smooth Deligne-Mumford stack.  We'll only deal with quasismooth hypersurfaces in this paper, and we'll furthermore assume that $C_f^*$ is not a {\it linear cone}, i.e., $f$ does not contain the monomial $x_i$ for any $i$.  This has the particular consequence that $C_f^*$ is not contained in any coordinate hyperplane of $\mathbb{A}^{n+2}$, and so $\mathcal{X} = [C_f^*/G]$ has trivial generic stabilizer.  It therefore makes sense to define the \textit{minimal log discrepancy} of $\mathcal{X}$ as the mld of the associated pair $(X,D)$, where $X \subset Y(\Sigma)$ is the coarse moduli space of $X$ and, as usual, $D$ has coefficient $1-\frac{1}{b}$ on the image of a prime divisor where the subgroup of $G$ acting trivially is of order $b$.

We'll use the following criterion for quasismoothness in linear systems, which is a generalization of a very similar statement due to Iano-Fletcher \cite[Theorem 8.1]{Iano-Fletcher}:

\begin{proposition}
\label{qsmoothcriterion}
Let $\mathcal{Y} = [(\mathbb{A}^{n+2} \setminus \{0\})/G]$ be a fake weighted projective stack and $T$ be a set of monomials which are homogeneous of degree $\chi$, where $\chi$ is not a coordinate character.  Then a general linear combination of monomials in $T$ defines a quasismooth hypersurface $\mathcal{X} \subset \mathcal{Y}$ if and only if for every nonempty set $I = \{i_1,\ldots,i_k \} \subset \{0,\ldots,n+1\}$, one of the following two conditions holds:
\begin{enumerate}
    \item[a.] there exists a monomial in $T$ of the form $x_{i_1}^{m_1} \cdots x_{i_k}^{m_k}$ (for some nonnegative $m_1,\ldots,m_k$), or
    \item[b.] there is another subset $J = \{j_1,\ldots,j_k\} \subset \{0,\ldots,n+1\}$ of size $k$ disjoint from $I$ such that for each $\mu = 1,\ldots,k$, there exists a monomial in $T$ of the form $x_{i_1}^{m_{1,\mu}}\cdots  x_{i_k}^{m_{k,\mu}} x_{j_{\mu}}$ (for some nonnegative $m_{1,\mu},\ldots,m_{k,\mu}$).
\end{enumerate}
\end{proposition}

\begin{proof}
For completeness, we'll include the proof, which proceeds along the same lines as Iano-Fletcher's.  Denote by $L$ the linear system of punctured affine cones $C_f^* \subset \mathbb{A}^{n+2} \setminus \{0\}$ of linear combinations $f$ of monomials in $T$.   By Bertini's theorem, a general member of $L$ is smooth away from the base locus $\mathrm{Bs}(L)$, which is a union of coordinate strata.  Therefore, it will suffice to check smoothness at points in the base locus. 

Given a stratum $W_I$ of dimension $k$, we may renumber indices so that $I = \{0,\ldots,k-1\}$.  If condition (a) holds for $I$, then there is some monomial in $T$ containing only the variables $x_0,\ldots,x_{k-1}$, so in particular a general $f$ is nonvanishing at any given point $p \in W_I$.  We've shown $W_I \cap \mathrm{Bs}(L) = \emptyset$ so a general $f$ intersects $W_I$ transversely and is smooth along this intersection.

In the event that (a) does not hold for $I$, we have $W_I \subset \mathrm{Bs}(L)$, but the condition (b) must hold.  Then we may expand $f$ as 

$$f = \sum_{i = k}^{n+1} x_i g_i(x_0,\ldots,x_{k-1}) + h,$$
where monomials in $h$ have total exponent of at least $2$ in the variables $x_k, \ldots, x_{n+1}$ (here we've used that there are no monomials in $T$ using only the first $k$ variables).  The partial derivatives of $f$ with respect to $x_0,\ldots,x_{k-1}$ vanish on $W_I$, but by condition (b), at least $k$ of the remaining partial derivatives $g_i$ are not identically zero on $W_I$.  Hence the locus in $W_I$ where the general $f$ is singular is the intersection of the base loci of $k$ free linear systems on $\mathbb{A}^{n+2}$ with $W_I$; this intersection has dimension at most $0$.  Since this locus is also $G$-invariant, it is empty.  Hence the general $C_f^*$ is smooth on $W_I$.

Conversely, if both conditions fail, then no linear combination of monomials in $T$ defines a quasismooth hypersurface.  Indeed, using the same expansion for $f$ above with respect to an $I$ for which (a) and (b) fail, there are fewer than $k$ nonvanishing $g_i$ in the sum for general $f$.  The intersection $Z \coloneqq \bigcap_{i = k}^{n+1} \{g_i = 0\} \cap W_I \subset W_I$ then has dimension at least $1$.  It follows that all derivatives of $f$ vanish on $Z$ for any linear combination $f$ of monomials in $T$, so $C_f^*$ is singular on the nonempty set $Z$.
\end{proof}

As above, let $T$ be a set of monomials of degree $\chi$ and let $L$ be the linear system spanned by $T$.  When the general member of $L$ is quasismooth, then we can describe the singularities of $\mathcal{X}$ very explicitly (cf. \cite[Proposition 2.6]{ETW}): 

\begin{proposition}
\label{fwphs_sing}
Let $\mathcal{Y} = [(\mathbb{A}^{n+2} \setminus \{0\})/G]$ be a fake weighted projective stack and $T$ be a set of monomials of degree $\chi$ spanning a linear system $L$.  Suppose that the general hypersurface $\mathcal{X}$ in $L$ (with coarse moduli space $X$) is quasismooth.   Let $p \in X$ be a closed point of the toric stratum $U_I \subset Y$ for $|I| = k$ and let $G_I = \bigcap_{i \in I} \ker(\theta_i)$.  Then,
\begin{enumerate}
    \item If $U_I$ is not in the base locus of $L$, then $\mathcal{X}$ has a quotient singularity of type $\frac{1}{G_I}(\theta_i|_{G_I}: i \notin I) \times \mathbb{A}^{k-2}$ at $p$.
    \item If $U_I$ is in the base locus of $L$, then there exists a $J \subset \{0,\ldots,n+1\}$ satisfying the conditions of \Cref{qsmoothcriterion}(b), and in particular, an index $j \in J$.  Then $\mathcal{X}$ has a quotient singularity of type $\frac{1}{G_I}(\theta_i|_{G_I}: i \notin I, i \neq j) \times \mathbb{A}^{k-1}$ at $p$.
\end{enumerate}
\end{proposition}

\begin{proof}
In both cases, the quasismoothness condition guarantees that $\mathcal{X}$ is locally given by a coordinate hyperplane slice through $p$ in $\mathcal{Y}$.  The resulting singularity will thus be the same as the singularity in $\mathcal{Y}$, with an appropriate weight removed.  Indeed, let $f$ be the polynomial defining $\mathcal{X}$.  In case (1), $C_f^*$ intersects the stratum $W_I$ transversely, so near a preimage of $p$ we can take it to have equation $x_i = 0$ for some $i \in I$.

In case (2), \Cref{qsmoothcriterion} guarantees that there is a set $J$ of $k$ indices so that the derivative $\frac{\partial f}{\partial x_j}$ doesn't vanish identically on $W_I$ for each $j \in J$.  There must be some $j \in J$ such that the corresponding partial derivative does not vanish at a preimage of $p$, so we can take the remaining coordinates as local coordinates of $C_f^*$ by the inverse function theorem.  Therefore, the quotient singularity of $\mathcal{X}$ at $p$ is the same as at $\mathcal{Y}$ with the $j$th coordinate character removed.  Furthermore, this singularity type is independent of which $j$ we choose: indeed, for every $j \in J$, by assumption $\theta_j^{-1} \chi$ is a product of nonnegative powers of the characters $\theta_i$ with $i \in I$.  Upon restricting to $G_I$, however, the latter characters become trivial.  Therefore, $\theta_j|_{G_I} = \chi|_{G_I}$ for any $j \in J$, as expected.
\end{proof}

\begin{remark}
\label{singularitiesremark}
As a particular consequence of this lemma, for any fixed toric stratum $U_I$, the hypersurface $\mathcal{X}$ has the same quotient singularity type at any intersection point with $U_I$.  Further, the worst singularities with respect to minimal log discrepancy occur on the smallest toric strata.  Therefore, to compute the mld of $\mathcal{X}$, it suffices to check singularities (1) on toric $1$-strata whose closures do not intersect $\mathrm{Bs}(L)$, and (2) at coordinate points ($0$-dimensional strata) in the base locus of $L$.
\end{remark}

\subsection{Berglund-H\"{u}bsch-Krawitz Mirror Symmetry}
\label{BHKintro}

Our main results on the minimal log discrepancies of certain Calabi--Yau pairs will be phrased in terms of mirror symmetry.  This takes advantage of a construction of mirror pairs due to Berglund-H\"{u}bsch-Krawitz (BHK) \cite{BH,Krawitz}.  We'll review the BHK mirror symmetry construction in this section, but see \cite{ABS} for further details.

Let $V \coloneqq \{f = 0\} \subset \mathbb{P}(a_0,\ldots,a_{n+1})$ be a hypersurface of dimension $n$ and weighted degree $d$, which is well-formed and quasismooth. Under these assumptions, $V$ is Calabi--Yau if and only if $d = a_0 + \cdots + a_{n+1}$.  Suppose that the weighted homogeneous polynomial (or \textit{potential}) $f$ defining the Calabi--Yau hypersurface $X$ has the same number of monomials as variables, namely $n+2$.  Then we may write
$$f = \sum_{i = 0}^{n+1} c_i \prod_{j = 0}^{n+1} x_j^{a_{ij}}.$$
The exponents $a_{ij}$ determine an $(n+2) \times (n+2)$ matrix $A$.  When this matrix is invertible, we say that $f$ is of \textit{Delsarte type} (this terminology will only be used for polynomials defining quasismooth Calabi-Yau hypersurfaces).  

When $f$ is of Delsarte type, it follows from \cite[Theorem 1]{KS} that $f$ can be written as a sum of \textit{atomic potentials} (up to coefficients).  There are three sorts of atomic potentials \cite[Section 2.2]{ABS}:
\begin{align*}
\label{atoms}
   & f_{\on{fermat}} = x^b, \\
   & f_{\on{loop}} = x_1^{b_1}x_2 + x_2^{b_2}x_3 + \cdots + x_{k-1}^{b_{k-1}}x_k + x_k^{b_k}x_1, \text{ and} \\
   & f_{\on{chain}} = x_1^{b_1}x_2 + x_2^{b_2}x_3 + \cdots + x_{k-1}^{b_{k-1}}x_k + x_k^{b_k}.
\end{align*}
The exponents $b_i$ in these equations are at least $2$, or else the degree would be the sum of just one or two weights, contradicting the fact that $V$ is Calabi--Yau. It's helpful to visualize potentials of Delsarte type as directed graphs, where an index $i$ points to $j$ if only if $x_i^{b_i}x_j$ is a monomial in $f$ (see \Cref{potentialgraph}).  This directed graph has the property that there is at most one arrow into and at most one arrow out of each node.  The connected components of the graph correspond to atomic potentials. In terms of the matrix $A$, the index $i$ points to an index $j$ in the graph of the potential if and only if $a_{ij} = 1$.

\begin{figure}
    \centering
    \begin{tikzcd}
& 5 \arrow[ld, bend right] & & & 3 \arrow[dd] & & \\
1 \arrow[dr, bend right] & & 7 \arrow[ul, bend right] & & & & 2 \\
& 4 \arrow[ur, bend right] & & & 6 & &
\end{tikzcd}
    \caption{The directed graph corresponding to a Delsarte potential function of shape $f = x_1^{b_1}x_4 + x_2^{b_2} + x_3^{b_3}x_6 + x_4^{b_4}x_7 + x_5^{b_5}x_1 + x_6^{b_6} + x_7^{b_7}x_5.$  This potential is composed of three atoms.} 
    \label{potentialgraph}
\end{figure}

From now on, suppose we're working with a $V$ defined by a potential of Delsarte type with associated matrix $A$. By the classification of atomic potentials, $A$ has nonnegative integer entries, the diagonal entries are at least $2$, and every row or column contains at most one nonzero off-diagonal entry, which must be a $1$. Without loss of generality, we can ignore the coefficients $c_i$ in $f$ and take them to be general nonzero constants.  This is because any two members of the linear system generated by the monomials $\prod_{j = 0}^{n+1} x_j^{a_{ij}}, i = 0,\ldots,n+1$, are isomorphic after multiplying by some element of the torus.   

The matrices corresponding to the atomic potentials $f_{\on{loop}}$ and $f_{\on{chain}}$ above are

$$A_{\on{loop}} = \begin{pmatrix}
b_1 & 1 & & &\\
& b_2 & 1 & &\\
& & \ddots & \ddots \\
 & & & b_{k-1} & 1 \\
1 & & & & b_k
\end{pmatrix} \text{ and }
A_{\on{chain}} = \begin{pmatrix}
b_1 & 1 & & &\\
& b_2 & 1 & &\\
& & \ddots & \ddots \\
 & & & b_{k-1} & 1 \\
 & & & & b_k
\end{pmatrix},$$
respectively.  We'll need the form of the inverses of these two types of matrices later, which may be readily computed via a cofactor expansion:

\begin{lemma}
\label{matrix_inverses}
    Let $A_{\on{loop}}$ and $A_{\on{chain}}$ be the matrices above.  Then

    $$A^{-1}_{\on{loop}} = \frac{1}{b_1 \cdots b_k + (-1)^{k-1}}\begin{pmatrix}
b_2 \cdots b_k & -b_3 \cdots b_k & \cdots & (-1)^{k-2} b_k & (-1)^{k-1}\\
(-1)^{k-1} & b_3 \cdots b_k b_1  & \cdots & (-1)^{k-3} b_k b_1 & (-1)^{k-2} b_1\\
\vdots & \vdots & \vdots & \vdots & \vdots \\
 b_2 \cdots b_{k-2} & \cdots & (-1)^{k-1}& b_k b_1 \cdots b_{k-2} & -b_1 \cdots b_{k-2}\\
-b_2 \cdots b_{k-1}  & \cdots & (-1)^{k-2} b_{k-1} & (-1)^{k-1} & b_1 \cdots b_{k-1}
\end{pmatrix}$$

$$\text{ and } A^{-1}_{\on{chain}} = \frac{1}{b_1 \cdots b_k}\begin{pmatrix}
b_2 \cdots b_k & -b_3 \cdots b_k & \cdots & (-1)^{k-2} b_k & (-1)^{k-1}\\
0 & b_3 \cdots b_k b_1  & \cdots & (-1)^{k-3} b_k b_1 & (-1)^{k-2} b_1\\
\vdots & \ddots & \ddots & \vdots & \vdots \\
0 & \cdots & 0 & b_k b_1 \cdots b_{k-2} & -b_1 \cdots b_{k-2}\\
0  & \cdots & 0 & 0 & b_1 \cdots b_{k-1}
\end{pmatrix}.$$
\end{lemma}

Define the \textit{charge} $q_i$ to be the sum of the entries of the $i$th row of $A^{-1}$.  Then the degree $d$ of $X$ is the least common denominator of the charges, and the weights satisfy $a_i = q_i d$ \cite[Section 2.2]{ABS}.  The transpose of $A$ defines a new potential $f^{\mathsf{T}}$.  The \textit{mirror charge} $q_i^{\mathsf{T}}$ is defined analogously as the sum of the entries of the $i$th column of $A^{-1}$.  Then $f^{\mathsf{T}}$ defines a Calabi--Yau hypersurface $X^{\mathsf{T}}$ with degree $d^{\mathsf{T}}$ equal to the least common denominator of the $q_i^{\mathsf{T}}$ in the weighted projective space $\mathbb{P}(a_0^{\mathsf{T}},\ldots,a_{n+1}^{\mathsf{T}})$, where $a_i^{\mathsf{T}} = d^{\mathsf{T}} q_i^{\mathsf{T}}$.  We'll refer to $d^{\mathsf{T}}$ and $a_0^{\mathsf{T}},\ldots,a_{n+1}^{\mathsf{T}}$ as the \textit{mirror degree} and \textit{mirror weights}, respectively.

Let $\mathrm{Aut}(f)$ be the group of diagonal automorphisms of $\C^{n+2}$ which preserve the potential $f$.  This is a finite group and is generated by the columns of $A^{-1}$, where a column $(c_0,\ldots,c_{n+1})^{\mathsf{T}}$ is interpreted as the diagonal automorphism $\mathrm{diag}(e^{2 \pi i c_0},\ldots,e^{2 \pi i c_{n+1}})$ {\cite[Section 3]{ABS}.  The action of the group $\mathrm{Aut}(f)$ descends to $X$ via the surjective homomorphism $\pi: \mathrm{Aut}(f) \rightarrow \mathrm{Aut}_T(V)$, where $\mathrm{Aut}_T(V)$ denotes the group of toric automorphisms of the weighted projective hypersurface $V$.  Let $J_f \coloneqq \ker(\pi)$ and $SL(f) \coloneqq \mathrm{Aut}(f) \cap \mathrm{SL}_{n+2}(\C)$.  To any group $J_f \subset G_f \subset SL(f)$, one can associate a group $G_f^{\mathsf{T}}$ satisfying $J_{f^{\mathsf{T}}} \subset G_f^{\mathsf{T}} \subset \mathrm{SL}(f^{\mathsf{T}})$.  The details of how to define this group won't be needed here.  Set $\widetilde{G_f} \coloneqq G_f/J_f \subset \mathrm{Aut}_T(V)$ and $\widetilde{G_f^{\mathsf{T}}} \coloneqq G_f^{\mathsf{T}}/J_{f^{\mathsf{T}}} \subset \mathrm{Aut}_T(V^{\mathsf{T}})$.  

We can now state the main theorem of BHK mirror symmetry (\cite[Theorem 2.3]{ABS}):

\begin{theorem}
The Calabi--Yau orbifolds $[V/\widetilde{G_f}]$ and $[V^{\mathsf{T}}/\widetilde{G_f^{\mathsf{T}}}]$ form a mirror pair, in the sense that 
$$H^{p,q}_{\on{CR}}([V/\widetilde{G_f}],\C) \cong H^{n-p,q}_{\on{CR}}([V^{\mathsf{T}}/\widetilde{G_f^{\mathsf{T}}}], \C),$$
where $H^{p,q}_{\on{CR}}(-,\C)$ denotes Chen-Ruan orbifold cohomology \cite{CR}.
\end{theorem}

\section{Minimal Log Discrepancies of Hypersurface Quotients}
\label{main_theorem_section}

In this section, we'll compute the minimal log discrepancies of certain klt Calabi--Yau pairs with standard coefficients which are quotients of hypersurfaces in weighted projective space.  It is conjectured that examples of this type achieve the smallest possible mld in each dimension, along with other extreme properties \cite[Conjectures 3.2, 3.4, 7.10]{ETWindex}; these conjectures are known to be true in low dimensions.  Our main result gives a connection between the mld's of these quotients and mirror symmetry:

\begin{theorem}
\label{mainmldthm}
Let $V_d \subset \mathbb{P}(a_0,\ldots,a_{n+1})$ be a well-formed quasismooth Calabi--Yau hypersurface defined by a polynomial of Delsarte type and $\mathrm{Aut}_T(V)$ the group of toric automorphisms of $V$.  Then, the minimal log discrepancy of the quotient pair $V/\mathrm{Aut}_T(V)$ is
$$\mathrm{mld}(V/\mathrm{Aut}_T(V)) = \frac{\min\{a_0^{\mathsf{T}},\ldots,a_{n+1}^{\mathsf{T}}\}}{d^{\mathsf{T}}}.$$
\end{theorem}

\begin{proof}
Let $f$ be the Delsarte-type equation defining the Calabi--Yau hypersurface $V$.  We may assume that $f$ has general coefficients, since any two quasismooth members of the linear system spanned by the monomials of $f$ differ by an automorphism of the torus. Let $L$ be the linear system generated by monomials in $f$. Following the notation of \Cref{BHKintro}, let $A$ be the $(n+2) \times (n+2)$ matrix determined by the exponents in $f$.  

Since $\mathrm{Aut}_T(V)$ is a finite subgroup of the torus, we may view $V/\mathrm{Aut}_T(V)$ as a hypersurface $\mathcal{X} \coloneqq \{f = 0\}$ inside the fake weighted projective stack $\mathcal{Y} \coloneqq [\mathbb{P}(a_0,\ldots,a_{n+1})/\mathrm{Aut}_T(V)] \cong [(\mathbb{A}^{n+2} \setminus \{0\})/G]$.  Here $G$ is the subgroup of the torus $(\C^*)^{n+2} \subset \mathbb{A}^{n+2}$ generated by the two smaller subgroups $\C^* = \{(t^{a_0},t^{a_1},\ldots,t^{a_{n+1}}): t \in \C^* \}$ and $\mathrm{Aut}(f)$.  We can now study the singularities of the hypersurface $\mathcal{X}$ using toric geometry and the results of \Cref{fwps}.  We'll use $(X,D)$ to denote the pair associated to the stack $\mathcal{X}$.

To understand the quotient morphism $\mathbb{A}^{n+2} \setminus \{0\} \rightarrow \mathcal{Y}$ of stacks from the standpoint of toric geometry, begin with the cone in the lattice $\Z^{n+2}$ generated by standard basis vectors $e_0,\ldots,e_{n+1}$.  The fan consisting of this cone and all its subcones is the fan of the toric variety $\mathbb{A}^{n+2}$.  The subvariety $\mathbb{A}^{n+2} \setminus \{0\} \subset \mathbb{A}^{n+2}$ corresponds to the subfan with the unique top-dimensional cone removed.  In \Cref{BHKintro}, we saw that the finite group $\mathrm{Aut}(f) \subset (\C^*)^{n+2}$ of automorphisms of the potential $f$ is generated by the columns $v_0, \ldots, v_{n+1}$ of $A^{-1}$.  Let $N$ be the superlattice of $\Z^{n+2}$ generated by the column vectors $v_0, \ldots, v_{n+1}$.  Then the datum for the fake weighted projective stack $\mathcal{Y} = \mathcal{Y}(\mathbf{\Sigma})$ is the triple $\mathbf{\Sigma} = (N', \Sigma, \{\bar{e}_0, \ldots, \bar{e}_{n+1}\})$, where $N' = N/(N \cap \mathrm{span}_{\R}\{(a_0,\ldots,a_{n+1})\})$, $\Sigma$ is the fan spanned by the image of the coordinate simplex $\mathrm{conv}(e_0,\ldots,e_{n+1})$ under the quotient $N_{\R} \rightarrow N'_{\R}$, and the $\bar{e}_i$ are the images of the basis vectors $e_i$.  The fan $\Sigma$ in $N'$ also defines the coarse moduli space $Y$ of $\mathcal{Y}$ as an ordinary toric variety.

Since $\mathcal{X}$ is a quasismooth hypersurface in $\mathcal{Y}$, the singularities of $\mathcal{X}$ will be suitable hyperplane slices of the quotient singularities of $\mathcal{Y}$.  By \Cref{fwphs_sing} and \Cref{singularitiesremark}, we only need to check the singularities of $\mathcal{X}$ on two types of toric strata to compute the mld of $\mathcal{X}$: the $1$-dimensional strata for which neither neighboring coordinate point is contained in the base locus of $L$, and the $0$-dimensional strata (coordinate points) which \textit{are} in the base locus of $L$.  We'll consider these two cases separately.  The key idea of the proof is that the images of the lattice points $v_i$ will always ``be responsible for" the smallest log discrepancy that occurs.

\textbf{Case 1}: Suppose that $\ell \cong \C^*$ is an open $1$-stratum of $Y$ whose closure does not intersect the base locus of $L$.  If the coordinates of the stratum are $x_i$ and $x_j$, this is the same as saying that in the graph of the corresponding Delsarte equation, neither $i$ nor $j$ ``points" at any other index.  For instance, in the example of \Cref{potentialgraph}, the $1$-stratum where only $x_2$ and $x_6$ are nonzero satisfies this condition.

Up to rearranging indices, we may assume that $\ell \subset Y$ is the locus where only the coordinates $x_n$ and $x_{n+1}$ are nonzero.  The cone in $\Sigma$ corresponding to this $1$-stratum is spanned by $\bar{e}_0, \ldots, \bar{e}_{n-1}$ and it is an $n$-dimensional cone in the $(n+1)$-dimensional vector space $N'_{\R}$.  We claim that $\mathrm{span}_{\R}\{\bar{v}_0,\ldots,\bar{v}_{n-1}\} = \mathrm{span}_{\R}\{\bar{e}_0,\ldots,\bar{e}_{n-1}\}$ and that $\mathrm{cone}\{ \bar{e}_0, \ldots, \bar{e}_{n-1} \} \subset \mathrm{cone} \{ \bar{v}_0, \ldots, \bar{v}_{n-1} \}$. (As the notation suggests, we're using $\bar{v}_i$ to mean the images of the columns $v_i$ of $A^{-1}$ in $N'$.) To prove this, we use that the vectors $\bar{e}_0, \ldots, \bar{e}_{n-1}$ are linearly independent in $N'_{\R}$ and each of $\bar{e}_0, \ldots, \bar{e}_{n-1}$ is a nonnegative $\Z$-linear combination of $\bar{v}_0, \ldots, \bar{v}_{n-1}$, via the following matrix equation:
\begin{equation}
\label{matrixlincombo}
   [e_0 \cdots e_{n+1} ] = I_{n+2} = A^{-1} A = [ v_0 \cdots v_{n+1} ] A.
\end{equation}
Expanding out the last product gives each standard basis vector as a nonnegative integer linear combination of columns $v_0,\ldots,v_{n+1}$.  Moreover, the coefficient of $v_i$ in the combination for $e_j$ is $a_{ij}$. The last two rows of $A$ have entries only on the diagonal by the assumption that neither $n$ nor $n+1$ points to any other index.  Therefore, $v_n$ and $v_{n+1}$ each have coefficient zero in the linear combination describing each of $e_0, \ldots, e_{n-1}$.  This completes the proof of the claim above.

We may add on either $\bar{v}_n$ or $\bar{v}_{n+1}$ to the set $\{\bar{v}_0, \ldots, \bar{v}_{n-1}\}$ to obtain a complete basis for the lattice $N'$ (either will suffice since $\bar{v}_0 + \cdots + \bar{v}_{n+1} = 0$).  Putting all this together, we've now described the singularities of $\mathcal{Y}$ on $\ell$ as $Q \times \mathbb{A}^1$, where $Q$ is the affine quotient singularity defined by the $\mathrm{cone}\{\bar{e}_0, \ldots, \bar{e}_{n-1}\}$ in the sublattice $K \subset N'$ with basis $\bar{v}_0, \ldots, \bar{v}_{n-1}$ and the extra $\mathbb{A}^1$ direction is given by the inclusion $K \subset N'$.  By \Cref{fwphs_sing}(1), the singularity of $\mathcal{X}$ at any intersection point with $\ell$ is simply $Q$.

The log discrepancy function for the affine singularity $Q$ is a linear function on $N'_{\R}$ which must be equal to $1$ at each of the distinguished points $\bar{e}_0, \ldots, \bar{e}_{n-1}$ on the respective rays.  This means that it lifts to some linear functional on the original space $\mathbb{R}^{n+2}$ with value $1$ at each of $e_0,\ldots,e_{n-1}$. So if $p \in \mathbb{R}^{n+2}$ is any (real) linear combination of $e_0, \ldots, e_{n-1}$, the log discrepancy function evaluated at the image of $p$ in $N'_{\R}$ equals the sum of coordinates of $p$ in $\R^{n+2}$ with the standard basis. By our work above, this applies in particular to any linear combination of $v_0, \ldots, v_{n-1}$.  Since the only lattice points in $N'$ inside $\mathrm{cone}\{\bar{e}_0, \ldots, \bar{e}_{n-1}\}$ are nonnegative integer combinations of $\bar{v}_0, \ldots, \bar{v}_{n-1}$, the log discrepancy of any divisor on a log resolution of $Q$ must be a nonnegative linear combination of the sums of coordinates of the $v_i$.  But sums of columns of $A^{-1}$ are simply the mirror charges $q_i^{\mathsf{T}}$.  We conclude that the log discrepancy of any divisor in a log resolution of the singularity $Q$ equals some nonnegative linear combination of the mirror charges.

\textbf{Case 2}: Suppose that $p$ is a $0$-dimensional stratum of $Y$ contained in the base locus of $L$.  This means that in the graph of the Delsarte equation, the corresponding index points to some other one.  For instance, in the example of \Cref{potentialgraph}, the coordinate points of $x_1, x_3, x_4, x_5,$ and $x_7$ all satisfy this condition.

Up to rearranging indices, we'll suppose that $p$ is the coordinate point of $x_n$ and that $n$ points to $n+1$.  This is the same thing as saying there is a monomial of the form $x_{n}^{b_{n}}x_{n+1}$ in $f$.  By the classification of Delsarte potentials, $n$ can only point at $n+1$, so $\mathcal{X}$ is defined in an \'{e}tale neighborhood of $p$ by the equation $x_{n+1} = 0$.  In the language of toric geometry, the affine quotient singularity $Q$ of $\mathcal{Y}$ at $p$ is the toric variety associated to the top-dimensional cone $\mathrm{cone}\{\bar{e}_0, \ldots, \bar{e}_{n-1}, \bar{e}_{n+1}\}$ in the fan of $\mathcal{Y}$.  By the discussion above (and in agreement with \Cref{fwphs_sing}(2)), the singularity of $\mathcal{X}$ at $p$ is the hyperplane slice $x_{n+1} = 0$ of $Q$.  This slice is another affine quotient singularity $Q'$, which as a toric Deligne-Mumford stack corresponds to the datum $(P,\mathrm{cone}\{f_0, \ldots, f_{n-1}\}, \{f_0, \ldots, f_{n-1} \})$.  To avoid cumbersome notation, we've used $P$ to represent the quotient lattice $N/(N \cap \mathrm{span}_{\R}\{(a_0, \ldots, a_{n+1}), e_{n+1}\})$. (Here we've taken the quotient by $e_{n+1}$ to represent the hyperplane slice $x_{n+1} = 0$.)  The point $f_i$ is the image in $P$ of $e_i$.  Similarly, we'll let $w_i$ be the image of $v_i$ in $P$.

Since $v_0 + \cdots + v_{n+1} = (q_0, \ldots, q_{n+1}) \in \mathrm{span}_{\R} \{a_0, \ldots, a_{n+1}\}$ and $v_n + b_{n+1} v_{n+1} = e_{n+1}$ (using the existence of the monomial $x_{n}^{b_{n}}x_{n+1}$ in $f$ and \eqref{matrixlincombo}), we have after the quotient that $w_0 + \cdots + w_{n+1} = w_n + b_{n+1} w_{n+1} = 0$.  Solving these equations gives that $w_n$ is a $\Z$-linear combination of the remaining vectors and
$$w_{n+1} = \frac{1}{b_{n+1}-1}(w_0 + \cdots + w_{n-1}).$$
Therefore, $w_0, \ldots, w_{n-1}$ span $P_{\R}$ as a vector space, but we might also have to add $w_{n+1}$ to get a full generating set for the lattice $P$ (if $b_{n+1} \geq 3)$.  In any case, the equation above shows $w_{n+1} \in \mathrm{cone}\{w_0, \ldots, w_{n-1}\}$.  We also claim that $\mathrm{cone}\{f_0, \ldots, f_{n-1}\} \subset \mathrm{cone}\{w_0, \ldots, w_{n-1}\}$.  Indeed, using \eqref{matrixlincombo} again, each $e_i$ is a linear combination of at most two of the columns of $A^{-1}$, and the only basis vectors which require $v_n$ in the combination are $e_n$ and $e_{n+1}$.  Thus, $f_0, \ldots, f_{n-1}$ are nonnegative $\Z$-linear combinations of $w_0, \ldots, w_{n-1}, w_{n+1}$, and we've already shown that $w_{n+1}$ is in the cone generated by the others.

It follows that any lattice point $x$ in $P \cap \mathrm{cone}\{f_0, \ldots, f_{n-1}\}$ is expressible as a nonnegative $\Z$-linear combination of $w_0, \ldots, w_{n-1}, w_{n+1}$. To see this, let $x = \gamma_0 w_0 + \ldots + \gamma_{n-1} w_{n-1}$ be the expression for $x$ in terms of the $\R$-basis $\{w_0, \ldots, w_{n-1}\}$.  Because of the inclusion $\mathrm{cone}\{f_0, \ldots, f_{n-1}\} \subset \mathrm{cone}\{w_0, \ldots, w_{n-1}\}$, all $\gamma_i$ are nonnegative.  If they are not integers, they all have identical fractional part which is some multiple of $\frac{1}{b_{n+1}-1}$, say $\frac{s}{b_{n+1}-1}$, because of our formula for $w_{n+1}$.  Thus,
$$x = \lfloor \gamma_0 \rfloor w_0 + \cdots \lfloor \gamma_{n-1} \rfloor w_{n-1} + s w_{n+1}$$
expresses $x$ as a nonnegative $\Z$-linear combination, as required.  We'll need the following fact about the special point $w_{n+1}$ later:

\begin{lemma}
\label{pointincone}
    $w_{n+1} \in \mathrm{cone}\{f_0, \ldots, f_{n-1}\}$.
\end{lemma}

\begin{proof}
In order to prove that $w_{n+1}$ is in $\mathrm{cone}\{f_0, \ldots, f_{n-1}\}$, we'll show that the coordinate functional associated to each coordinate $f_i$ in the basis $\{f_0,\ldots,f_{n-1}\}$ of $P_{\R}$ is positive at $w_{n+1}$.  It's easier to view these as induced by linear functionals on the original space $\R^{n+2}$ of which $P_{\R}$ is a quotient.

Indeed, for each $i = 0,\ldots,n-1$, let $\pi_i: \R^{n+2} \rightarrow \R$ be the unique linear functional descending to the coordinate functional for $f_i$ on $P_{\R}$.  As above, $e_0, \ldots, e_{n+1}$ are the standard basis vectors for $\R^{n+2}$.  Then $\pi_i$ is determined by the conditions: $\pi_i(e_i) = 1$, $\pi_i(e_j) = 0$ for $i \neq j$ and $i \neq n$, and $\pi(a) = 0$, where $a = (a_0, \ldots, a_{n+1})$ is the vector of weights.  We can use these conditions to determine the value of $\pi_i$ on the final basis vector $e_n$:
$$\pi_i(e_n) = -\pi_i\left(\frac{1}{a_n}(a_0e_0 + \cdots + a_{n-1} e_{n-1} + a_{n+1} e_{n+1})\right) = -\frac{a_i}{a_n}.$$
We can conclude that for a vector $r = (r_0, \ldots,r_{n+1}) \in \R^{n+2}$, $\pi_i(r) \geq 0$ if and only if $r_i - \frac{a_i}{a_n}r_n \geq 0$, or in other words, $r_i \geq \frac{a_i}{a_n} r_n$.  To complete the proof of the lemma, we need to verify this inequality for $r = v_{n+1}$ specifically, since this vector maps to $w_{n+1}$ under the quotient $\R^{n+2} \rightarrow P_{\R}$.  

We know that $n$ points to $n+1$ and that both are part of an atom which is either a chain or a loop potential.  The entry $r_n$ of $v_{n+1}$ that we're interested in corresponds to a superdiagonal element of one of the matrix inverses in \Cref{matrix_inverses}; therefore, $r_n$ is negative.  It follows that if the index $i$ is part of a distinct atom from $n$ and $n+1$ (so that $r_i = 0$), the required inequality is automatically satisfied.  We can then reduce to the case of looking at only loop or chain atoms.

The required inequality in these two special cases follows directly from \Cref{matrix_inverses}; we'll switch to using the notation of that lemma to discuss the various column entries.  In the loop case, this means that we'll look at the last column $v_k = (r_1,\ldots,r_k)$ of the matrix $A^{-1}_{\on{loop}}$ in \Cref{matrix_inverses} and try to prove $r_i \geq \frac{a_i}{a_{k-1}} r_{k-1}$ for $i = 1,\ldots,k-2$. The easiest way to see this is to observe that, modulo $\Z$, the $i$th row of $A^{-1}_{\on{loop}}$ is simply the previous row times $-b_{i-1}$ for each $i$.  Therefore, the same fact holds true for the entries of the vector which is the sum of the columns, namely the charge vector $q = (a_1/d,\ldots,a_k/d)$.  We can observe that the $(k-1)$-entry $r_{k-1}$ of the $k$th column of $A^{-1}_{\on{loop}}$, namely $-b_1 \cdots b_{k-2}$, is $(-1)^{k-i-1} b_i \cdots b_{k-2}$ times the $i$th entry $r_i$ of the same column for all $i = 1,\ldots,k-2$.  Since the charge vector has a similar property modulo $\Z$, $(-1)^{k-i-1} b_i \cdots b_{k-2} a_i \equiv a_{k-1} \pmod d$.  If $(-1)^{k-i-1} b_i \cdots b_{k-2}$ is negative, the $i$th entry is positive and the required inequality will hold automatically, while if $(-1)^{k-i-1} b_i \cdots b_{k-2}$ is positive, $b_i \cdots b_{k-2}$ is at least as large as $a_{k-1}/a_i$, since all weights are between $0$ and $d$ and $a_i b_i \cdots b_{k-2} \equiv a_{k-1} \pmod d$ by the discussion above.  This finishes the proof of the required inequalities. 

The argument for the chain potential case is very similar, because $A^{-1}_{\on{chain}}$ has nearly the same properties as $A^{-1}_{\on{loop}}$.  Indeed, above the main diagonal, consecutive rows have the same relation as in $A^{-1}_{\on{loop}}$, and the same argument works. Below the main diagonal, entries of the column of interest might be zero, but the inequality again holds automatically in this case. \end{proof}

Finally, we need to compute the log discrepancy function for the singularity $Q'$:

\begin{lemma}
\label{logdiscrepfunction}
    The log discrepancy function $\psi$ for $Q'$ has value $q_i^{\mathsf{T}}$ on $w_i$ for $i = 0, \ldots, n-1, n+1$.
\end{lemma}

\begin{proof}
The singularity $Q'$ is affine, so the log discrepancy function $\psi$ is simply the unique linear functional on $P_{\R}$ with the property that it has value $1$ on each of the vectors $f_0, \ldots, f_{n-1}$.  We'll proceed in much the same way as the previous lemma by doing our calculations on $\R^{n+2}$; the linear functional $\psi$ on the quotient vector space $P_{\R}$ is induced by a unique linear functional $\tilde{\psi}$ on $\R^{n+2}$ with the property that $\tilde{\psi}(e_0) = \cdots = \tilde{\psi}(e_{n-1}) = 1$, $\tilde{\psi}(a) = 0$, and $\tilde{\psi}(e_{n+1}) = 0$, where $a = (a_0, \ldots, a_{n+1})$ is the vector of weights.  We're looking to compare this to the ``sum of coordinates" linear functional $\phi$, which has a value of $1$ at each of $e_0, \ldots, e_{n+1}$.  First, using $\tilde{\psi}(a) = 0$, we compute
$$\tilde{\psi}(e_{n}) = -\tilde{\psi}\left(\frac{1}{a_n}(a_0e_0 + \cdots + a_{n-1} e_{n-1} + a_{n+1} e_{n+1})\right) = - \frac{a_0 + \cdots + a_{n-1}}{a_n}.$$
For any point $r = (r_0,\ldots,r_{n+1}) \in \R^{n+2}$, the linear functionals $\tilde{\psi}$ and $\phi$ coincide on $r$ if and only if
$$r_0 + \cdots + r_{n+1} = r_0 + \cdots + r_{n-1} - \frac{a_0 + \cdots + a_{n-1}}{a_n}r_n.$$
Simplifying this gives
$$r_{n+1} = - \frac{a_0 + \cdots + a_{n-1} + a_n}{a_n}r_n = - \frac{d-a_{n+1}}{a_n} r_n = -b_n r_n,$$
where we've used that $a_0 + \cdots + a_{n+1} = d$ by the Calabi--Yau property of $V$ and that $x_n^{b_n} x_{n+1}$ is a monomial of weighted degree $d$ with these weights.  

We'll show that $\phi = \tilde{\psi}$ holds for $v_0, \ldots, v_{n-1}$, and $v_{n+1}$ (and hence also for any linear combinations of these).  The criterion $r_{n+1} = -b_n r_n$ only involves the last two coordinates of $r$, so it is automatically satisfied for any $v_i$ with $i$ not in the same atomic potential as $n$ or $n+1$ (since the relevant coordinates in the block diagonal matrix $A^{-1}$ will both be zero in that case).

Thus, it suffices to verify this condition for the relevant columns of the inverse matrix $A^{-1}$ where $A$ is the matrix of an atomic loop or chain potential.
 Indeed, one can readily check from \Cref{matrix_inverses} that the $i+1,j$ entry of $A^{-1}$ is $-b_i$ times the $i,j$ entry, unless $i = j$.  This proves that $r_{n+1} = -b_n r_n$ holds for every column vector except $v_n$ itself.  Thus, for $i = 0, \ldots, n-1,n+1$, the value of $\tilde{\psi}$ on $v_i$ must equal the value of $\phi$ on $v_i$, namely the mirror charge $q_i^{\mathsf{T}}$. This means $\psi(w_i) = \tilde{\psi}(v_i) = q_i^{\mathsf{T}}$ for the required values of $i$, completing the proof.
\end{proof}

Since any point in $P \cap \mathrm{cone}\{f_0, \ldots, f_{n-1}\}$ is a nonnegative $\Z$-linear combination of $w_0, \ldots, w_{n-1}, w_{n+1}$, \Cref{logdiscrepfunction} shows that the log discrepancy of any divisor in a log resolution of the singularity $Q'$ is a nonnegative linear combination of the mirror charges.

In summary, the analysis in Cases (1) and (2) has shown that any log discrepancy of a divisor in a log resolution of the pair $(X,D)$ must be a nonnegative linear combination of mirror charges $q_i^{\mathsf{T}}$.  Since $(X,D)$ is the pair associated to $V/\mathrm{Aut}_T(V)$, we've shown
$$\mathrm{mld}(V/\mathrm{Aut}_T(V)) \geq \min\{q_0^{\mathsf{T}},\ldots,q_{n+1}^{\mathsf{T}}\}.$$
To demonstrate equality, the last step is to show that the smallest mirror charge actually appears as a log discrepancy of some divisor over $(X,D)$.  This will follow quickly from the work we've done already.  Indeed, suppose without loss of generality that $q_{n+1}^{\mathsf{T}}$ is the smallest mirror charge.  If any index (assume it's $n$) points to $n+1$, we saw in \Cref{pointincone} that the image $w_{n+1}$ of $v_{n+1}$ is a lattice point in the cone of the singularity of $\mathcal{X}$ at the coordinate point of $x_n$.  Therefore, we get a divisor of log discrepancy $q_{n+1}^{\mathsf{T}}$ by \Cref{logdiscrepfunction}.  Otherwise, no other index points to $n+1$.  Representing $A$ as a block diagonal matrix separated into atomic potentials, this means that the $n+1$ column of $A$ has only the diagonal entry $b$, which is the top left corner of an upper triangular block inside $A$.  It follows that $v_{n+1} = \frac{1}{b} e_{n+1}$ and the subgroup of $\mathrm{Aut}_T(V)$ acting as the identity on the divisor $\{x_{n+1} = 0\}$ has order $b$.  In the coarse moduli pair $(X,D)$, the divisor $\{x_{n+1} = 0\} \subset X$ must have coefficient $1 - \frac{1}{b}$ in $D$.  This divisor therefore has log discrepancy of 
$$1 - \left(1 - \frac{1}{b} \right) = \frac{1}{b} = q_{n+1}^{\mathsf{T}},$$
as required.
\end{proof}

\section{Applications}
\label{applications_section}

In this section, we use \Cref{mainmldthm} to construct Calabi--Yau varieties or pairs with particularly extreme properties. We'll first examine how the main theorem offers an alternative perspective on an example due to Jihao Liu \cite[Remark 2.6]{Liu} of a klt Calabi--Yau pair with standard coefficients and small mld.  Later in the section, we'll construct examples of klt Calabi--Yau varieties (rather than pairs) with conjecturally minimal mld.

Liu's example is closely related to another by Koll\'{a}r of a klt pair with standard coefficients such that $K_X + D$ is ample and has very small volume (see \cite{Kollarlog}, \cite[Introduction]{HMXbir}).  The Calabi--Yau pair is defined as follows: 
$$(X,D) = \left(\mathbb{P}^n, \frac{1}{2}H_0 + \frac{2}{3}H_1 + \cdots + \frac{s_n - 1}{s_n}H_n + \frac{s_{n+1} - 2}{s_{n+1} - 1}H_{n+1} \right).$$
Here $H_0, \ldots, H_{n+1}$ are general hyperplanes in $\mathbb{P}^n$, and the $s_i$ denote \textit{Sylvester's sequence}.  This sequence begins $s_0 = 2$ and is defined by $s_m = s_0 \cdots s_{m-1} + 1$ for all $m \geq 1$.  The sequence grows doubly exponentially and its terms are pairwise relatively prime; the first few terms are $2,3,7,43,1807,\ldots$.  It's clear that $(X,D)$ is klt, since it is its own log resolution.  It also has standard coefficients, $K_X + D \sim_{\Q} 0$, and the smallest log discrepancy comes from the divisor $H_{n+1}$, so $\mathrm{mld}(X,D) = 1/(s_{n+1}-1)$.  It is conjectured \cite[Conjecture 3.2]{ETWindex} that this is the smallest possible mld for any Calabi--Yau klt pair with standard coefficients and dimension $n$. This conjecture is proven in dimensions at most $2$ (see \cite[Proposition 6.8]{ETWindex}).

The pair $(X,D)$ above has another interpretation as a quotient of the hypersurface $V \subset \mathbb{P}^{n+1}(d/s_0,d/s_1,\ldots,d/s_n,1)$ of degree $d = s_{n+1} - 1$ defined by the Fermat-type potential equation
$$f \coloneqq x_0^2 + x_1^3 + \cdots + x_n^{s_n} + x_{n+1}^{s_{n+1}-1} = 0.$$
The automorphism group of the potential $f$ is isomorphic to $\mu_2 \oplus \mu_3 \oplus \cdots \oplus \mu_{s_n} \oplus \mu_{s_{n+1}-1}$, where the $i$th summand acts on $x_i$ by multiplication by the corresponding root of unity and fixes all other variables.  The toric automorphism group $\mathrm{Aut}_T(V)$ is the quotient of this group by the action of the subgroup $\C^* \cap \mathrm{Aut}(f) \cong \mu_{s_{n+1}-1}$.  We have $\mathrm{Aut}_T(V) \cong \mu_{s_{n+1}-1}$, generated by 
$$\zeta \cdot (x_0: x_1: \cdots :x_n: x_{n+1}) = (\zeta^{d/s_0}x_0: \zeta^{d/s_1}x_1: \cdots : \zeta^{d/s_n}x_n: x_{n+1}),$$
where $\zeta \in \mu_{s_{n+1}-1}$.  There is an isomorphism of varieties $V/\mathrm{Aut}_T(V) \cong \mathbb{P}^n$ given by 
$$(x_0:x_1: \cdots :x_n:x_{n+1}) \mapsto (x_0^{2}: x_1^{3}: \cdots : x_n^{s_n}:x_{n+1}^{s_{n+1}-1}),$$
which identifies $V/\mathrm{Aut}_T(V)$ with the hyperplane $y_0 + y_1 + \cdots + y_n + y_{n+1} = 0$ in $\mathbb{P}^{n+1}$.  This hyperplane is, of course, isomorphic to $\mathbb{P}^n$.  The pair structure on the quotient is precisely Liu's pair above.

\Cref{mainmldthm} reproduces the same value for the mld of this pair.  Indeed, the mirror weights and degree of $V$ equal the original weights and degree since $f$ is Fermat (and so the corresponding matrix $A$ from \Cref{BHKintro} is diagonal).  The smallest mirror weight is therefore $1$, so the mld of $V/\mathrm{Aut}_T(V)$ is $1/(s_{n+1}-1)$, as expected.

Next, we construct klt Calabi--Yau varieties (rather than pairs) with very small mld.  The construction shares some similarities with the large index Calabi--Yau example due to Totaro, Wang, and the author \cite[Section 7]{ETWindex}.  The intricate identities relating the constants defining both of these examples are worked out fully in \cite{ET}. Therefore, we only define the examples below and sketch proofs of their properties.  The key idea is to begin with a loop potential (resp. a potential of the form $x_0^2 + [\on{loop}]$), so that the action of the toric automorphism group (resp. an index $2$ subgroup of the toric automorphism group) acts freely in codimension $1$.  After taking the quotient, we then obtain a variety, rather than a pair.  The example will be defined recursively using Sylvester's sequence. The formulas are rather complicated, but we write down explicit examples for $n = 2,3,4$ at the end of this section.

For brevity, given integers $b_{i_1}, \ldots, b_{i_k}$, we introduce the symbol $B_{i_1, \ldots, i_k}$ for the alternating sum $b_{i_1}b_{i_2} \cdots b_{i_k} - b_{i_1} \cdots b_{i_{k-1}} + \cdots + (-1)^{k-1} b_{i_1} + (-1)^k$.

\begin{definition}
\label{kltexampledef1}
For $n = 2r + 1, r \geq 1$ or $n = 2r$, $r \geq 1$, we'll define integers $b_0, \ldots, b_n, v_n$ as follows.  When $0 \leq i \leq r$, define $b_i \coloneqq s_i$.  Then, define the remaining $b_i$ inductively via the formulas
\begin{equation*}
b_{r+i} \coloneqq 1 + (b_{r+1-i}-1)^2 B_{r+1,r,r+2,r-1,\ldots,r+i-1,r-i+2}
\end{equation*}
for $1\leq i\leq r+1$ when $n=2r+1$ or $1\leq i\leq r$ when $n=2r$.

For $n = 2r + 1$, $r \geq 1$, define an integer $v_{2r+1}$ by:
$$v_{2r+1} \coloneqq B_{r+1,r,r+2,r-1,\ldots,2r+1,0}-B_{r,r+2,r-1,\ldots,2r+1,0}+\cdots-B_0.$$
For $n = 2r$, $r \geq 1$, define an integer $v_{2r}$ by:
$$v_{2r} \coloneqq 2(B_{r+1,r,r+2,r-1,\ldots,2r,1}-B_{r,r+2,r-1,\ldots,2r,1}+\cdots-B_1)+1.$$
\end{definition}

Note that a $B$ with empty subscript is considered to be $1$, so $b_{r+1} = 1 + (b_r-1)^2$. These formulas define the exponents of the equation for our hypersurface example.  The equations have the following form. When $n=2r+1$ for $r\geq 1$,
let the hypersurface $V$ have equation
\begin{equation}\label{oddeq}
x_0^{b_0}x_{2r+2}+ x_1^{b_1} x_{2r+1} + \cdots + x_r^{b_r}x_{r+2}+x_{r+1}^{b_{r+1}}x_r + x_{r+2}^{b_{r+2}} x_{r-1} + \cdots +x_{2r+1}^{b_{2r+1}}x_0
       +x_{2r+2}^{v_{2r+1}}x_{r+1} = 0.
\end{equation}
The potential defining this hypersurface is atomic of loop type. Similarly, when $n = 2r$ for $r \geq 1$,
let $V$ have equation given by
\begin{equation}\label{eveneq}
x_0^{b_0} + x_1^{b_1}x_{2r+1}+ x_2^{b_2} x_{2r} + \cdots + x_r^{b_r}x_{r+2}+x_{r+1}^{b_{r+1}}x_r
       + x_{r+2}^{b_{r+2}} x_{r-1} + \cdots +x_{2r}^{b_{2r}}x_1+x_{2r+1}^{v_{2r}}x_{r+1} = 0.
\end{equation}
The potential in this case is of the form $x_0^2 + [\on{loop}]$.

The weights $a_0,\ldots,a_{n+1}$ and degree $D$ of the hypersurface $V \subset \mathbb{P}^{n+1}(a_0,\ldots,a_{n+1})$
are uniquely determined by the equations above, once we require that $\gcd(a_0,\ldots,a_{n+1}) = 1$, which is necessary for the weighted projective space to be well-formed.  The form of the equations guarantees that \eqref{oddeq} and \eqref{eveneq} define quasismooth hypersurfaces of each dimension $n$. If $n = 2r + 1$ is odd, the toric automorphism group $G \coloneqq \mathrm{Aut}_T(V)$ acts freely in codimension $1$ on $V$ because it is defined by a loop potential \cite[Proposition 7.2]{ETWindex}.  Therefore, the quotient $V/G$ is a variety, rather than a pair.

When $n = 2r$ is even, the hypersurface $V$ is not defined by a loop potential. Instead of taking the quotient of $V$ by $\mathrm{Aut}_T(V)$, we instead consider $V/G$, where $G$ is now the index $2$ subgroup of $\mathrm{Aut}_T(V)$ corresponding to the loop atom.  Using the same reasoning as in \cite[Proposition 7.2]{ETWindex}, one can check that $G$ acts freely in codimension $1$ on $V$. On its own, \Cref{mainmldthm} does not say anything about $V/G$ since $G$ is not the full toric automorphism group of $V$.  However, given that $V/\mathrm{Aut}_T(V) \cong (V/G)/\mu_2$, one can deduce in this case that the mld of $V/G$ is twice that of its quotient.

We define the additional constants:

\begin{definition}
\label{kltvarietydef2}
$$D \coloneqq \begin{cases}
D_{2r+1} = B_{r+1,r,r+2,r-1,\ldots,2r+1,0}&\text{if }n=2r+1,\\
D_{2r} = 2B_{r+1,r,r+2,r-1,\ldots,2r,1} &\text{if }n=2r.
\end{cases}$$ 

$$m \coloneqq \begin{cases} m_{2r+1} = B_{0,2r+1,1,2r,\ldots,r,r+1} &\text{if }n=2r+1,\\
m_{2r} = B_{1,2r,2,2r-1,\ldots,r,r+1} &\text{if }n=2r.
\end{cases}$$
\end{definition}

The properties of the klt variety example are summarized in the following theorem.

\begin{theorem}{\cite[Theorem 5.1]{ET}}
\label{exampletheorem}
In each dimension $n\geq 2$, the hypersurface $V$ defined above is 
well-formed, Calabi-Yau and quasismooth of degree $D$.
The hypersurface $V$ carries an action of the cyclic group $\mu_m$ such that $V/\mu_m$
is a complex klt Calabi--Yau variety
with mld $1/m$. This mld is smaller than $1/2^{2^n}$ for each $n > 2$.
\end{theorem}

The verifications of these properties are carried out fully in \cite{ET}, but the 
key point is showing the identities \cite[Proposition 4.1]{ET}
\begin{align}
\label{identities}
\begin{split}
   m_{2r+1}D_{2r+1}-1
& = b_0\cdots b_{2r+1}v_{2r+1},\\
m_{2r}D_{2r}-1
& = b_1\cdots b_{2r}v_{2r}. 
\end{split}
\end{align}

For instance, when $n = 2r+1$, the determinant of the loop matrix associated to the
potential equation \eqref{oddeq} is $b_0\cdots b_{2r+1}v_{2r+1} + 1$.
One can use \eqref{identities} to prove the degree of $V$ is $D$, the smallest weight
is $a_{n+1} = 1$, and the smallest mirror charge is equal to $1/m$.  The value of the mld is then a 
consequence of \Cref{mainmldthm}.  The situation for even $n$ is similar.

The value $1/m$ of the mld decays doubly exponentially with dimension and is extremely close to that of the conjecturally optimal klt pair mentioned earlier in the section.  Indeed, the mld value for varieties is within a constant factor of less than $6$ of the value for pairs when $n$ is even, and less than $23$ when $n$ is odd \cite[Lemma 8.1]{ET}. This is compelling evidence for the conjecture that this example is optimal for klt varieties:

\begin{conjecture}
\label{kltvarconj}
For every $n \geq 2$, the quotient $V/G$ defined above has the smallest mld of any klt Calabi--Yau variety of dimension $n$.
\end{conjecture}

This conjecture is also supported by additional evidence in small dimensions.  Indeed, in dimension $2$, $V$ is the degree $22$ hypersurface
$$\{x_0^2 + x_1^3 x_3 + x_1 x_2^5 + x_2 x_3^{19} = 0\} \subset \mathbb{P}^3(11,7,3,1).$$
The group $G$ of order $13$ acts freely in codimension $1$ on $V$, and the quotient variety $V/G$ has mld $\frac{1}{13}$.  This is the smallest possible for any klt Calabi--Yau surface by \cite[Proposition 6.1]{ETWindex}.  Note that the mld of the quotient $V/\mathrm{Aut}_T(V)$ by the full toric automorphism group, of order $26$, has mld $\frac{1}{26}$, as predicted by \Cref{mainmldthm}, since the mirror hypersurface to $V$ has degree $26$ in $\mathbb{P}(13,7,5,1)$.  However, that quotient is a pair, rather than a variety.  

In dimension $3$, $V$ is the degree $191$ hypersurface
$$\{x_0^2 x_4 + x_1^3 x_3 + x_1 x_2^5 + x_0 x_3^{12} + x_2 x_4^{165} = 0\} \subset \mathbb{P}^4(95,61,26,8,1).$$
The group $G$ of order $311$ acts freely in codimension $1$ on $V$, and \Cref{mainmldthm} shows that the quotient variety $V/G$ is a klt Calabi--Yau 3-fold with mld $\frac{1}{311}$.  In dimension $4$, $V$ is the degree $925594$ hypersurface
$$\{x_0^2 + x_1^3 x_5 + x_2^7 x_4 + x_2 x_3^{37} + x_1 x_4^{893} + x_3 x_5^{904149} = 0\} \subset \mathbb{P}^5(462797,308531,132129,21445,691,1).$$

The group $G$ of order $677785$ acts freely in codimension $1$ on $V$, and \Cref{mainmldthm} shows that the quotient variety $V/G$ is a klt Calabi--Yau 4-fold with mld $\frac{1}{677785}$.  For each of these computations, applying \Cref{mainmldthm} is a much simpler way of computing the mld than analyzing all the quotient singularities of $V/G$ individually.

In dimensions $3$ and $4$, our examples have the smallest mld of any possible examples of this type. That is, they are optimal among all klt Calabi--Yau varieties arising as quotients by toric automorphisms of quasismooth hypersurfaces in weighted projective space defined by Delsarte potentials.  This was verified by computer search, using the databases of Calabi--Yau 3-fold and 4-fold hypersurfaces in \cite{database}.

\end{document}